\documentclass[12pt,oneside,reqno]{article}
\usepackage{braket, amsxtra, xr}
\usepackage[utf8]{inputenc}

\usepackage{slashed}
\usepackage{graphicx}
\usepackage{xcolor}
\usepackage{pdfrender}
\pdfrender{TextRenderingMode=2,LineWidth=.5pt}
\usepackage{amsmath}
\usepackage{amsthm}
\usepackage{indentfirst}
\usepackage{amsfonts}
\usepackage{marginnote}
\usepackage{xcolor}
\usepackage{hyperref}
\usepackage{OldStandard}
\usepackage{multicol}
\usepackage{lipsum}
\usepackage{newtxmath}

\usepackage{geometry}

\usepackage{blindtext}
\usepackage{authblk}
\usepackage{mathtools}

\usepackage{xcolor}
\usepackage{pdfrender}
\pdfrender{TextRenderingMode=2,LineWidth=.5pt}

\usepackage{lipsum} 

\newtheoremstyle{plainindent}
  {\topsep}   
  {\topsep}   
  {\itshape}  
  {\parindent}
  {\bfseries} 
  {.}         
  {5pt plus 1pt minus 1pt} 
  {}          

\theoremstyle{plainindent}

\parskip=\baselineskip

\newcommand\myshade{85}
\colorlet{mylinkcolor}{blue}
\colorlet{mycitecolor}{red}
\hypersetup{linkcolor  = mylinkcolor,
	citecolor  = mycitecolor,
	urlcolor   = myurlcolor!\myshade!black,
	colorlinks = true,}
	
\title{\textsc{ functional equation for mellin transform of fourier
series associated with modular forms}}
\author{\textit{By} \textsc{ Omprakash Atale${}^{\dagger}$ }\footnote{${}^{\dagger}$E-mail: atale.om@outlook.com }}
\affil{\small ${}^{\dagger}$\textit{Department of Mathematics, Savitribai Phule Pune University,}\\ \textit{Pune-411007, India}}
\date{[10 September 2024]}
\begin{document}
\maketitle
\begin{abstract}
   Let $X_1(s)$ and $X_2(s)$ denote the Mellin transforms of $\chi_{1}(x)$ and $\chi_{2}(x)$, respectively. Ramanujan investigated the functions $\chi_1(x)$ and $\chi_2(x)$ that satisfy the functional equation
\begin{equation*}
X_{1}(s)X_2(1-s) = \lambda^2,
\end{equation*}
where $\lambda$ is a constant independent of $s$. Ramanujan concluded that elementary functions such as sine, cosine, and exponential functions, along with their reasonable combinations, are suitable candidates that satisfy this functional equation. Building upon this work, we explore the functions $\chi_1(x)$ and $\chi_2(x)$ whose Mellin transforms satisfy the more general functional equation
\begin{equation*}
\frac{X_1(s)}{X_2(k-s)} = \sigma^2,
\end{equation*}
where $k$ is an integer and $\sigma$ is a constant independent of $s$. 

As a consequence, we show that the Mellin transform of the Fourier series associated with certain Dirichlet L-functions and modular forms satisfy the same functional equation.
\end{abstract}
\newpage

\fontsize{11.5}{16}\selectfont
\begin{center}
    \textsc{\textbf{\large \S I. INTRODUCTION}}
\end{center}
Let 
\begin{equation*}
    X_1(s)=\int_{0}^{\infty}{x^{s-1}\chi_{1}(x)dx}\quad\mathrm{and}\quad X_2(s)=\int_{0}^{\infty}{x^{s-1}\chi_{2}(x)dx}.\tag{1.1}
\end{equation*}
Ramanujan [\cite{1}. pg. 223-227] considered the problem of identifying and constructing $\chi_{1}(x)$ and $\chi_{2}(x)$ such that the functional equation for the product of Mellin transforms 
\begin{equation*}
   X_1(s)X_2(1-s)=\lambda^2\tag{1.2} 
\end{equation*}
is satisfied. Here, $\lambda$ is a constant independent of $s$. Ramanujan's results were based on his claim that the two equations
\begin{equation*}
    \int_{0}^{\infty}\phi(x)\chi_1(nx)dx=\lambda \psi(n)\quad\mathrm{and}\quad \int_{0}^{\infty}\psi(x)\chi_2(nx)dx=\lambda \phi(n)\tag{1.3}
\end{equation*}
imply each other. The functions $\phi(x)$ and $\psi(x)$ are to be discussed later. Let $Z_1(s)$ and $Z_2(s)$ denote the Mellin transform of $\phi(x)$ and $\psi(s)$, respectively. Then Ramanujan showed that
\begin{equation*}
    \frac{Z_{1}(s)}{Z_2(1-s)}=\frac{X_2(s)}{\lambda}=\frac{\lambda}{X_1(1-s)}.\tag{1.4}
\end{equation*}

The functions $\chi_{1}(x)=\chi_2(x)=\sin(x)$ and $\chi_{1}(x)=\chi_2(x)= \cos(x)$ are some preliminary examples whose Mellin transform satisfies the functional equation for the product of Mellin transforms with $\lambda=\frac{\pi}{2}$. Other examples include the ordinary Bessel function $\chi_{1}(x)=\chi_2(x)=\sqrt{x}J_\nu(x)$ ($\Re(\nu)>-1$) with $\lambda=1$ and
\begin{equation*}
    \chi_{1}(x)=\chi_2(x)=\frac{x^\alpha}{1+x^2}\tag{1.5}
\end{equation*}
with $\lambda=\frac{\pi}{2}$ where $\alpha$ is an integer. Following is an example of Eqn. (1.3) 
\begin{equation*}
    \int_{0}^{\infty}e^{-x^{2}}\cos(2nx)dx=\frac{\sqrt{\pi}}{2} e^{-n^{2}}\quad\mathrm{and}\quad \int_{0}^{\infty}xe^{-x^{2}}\sin(2nx)dx=\frac{n\sqrt{\pi}}{2} e^{-n^{2}}\tag{1.6}.
\end{equation*}

Ramanujan also concluded that a reasonable combination of the functions $\sin(x)$, $\cos(x)$, and $e^{-x}$ can also be a solution to Eqn. (1.2). For example
\begin{equation*}
    \chi_{1}(x)=\chi_{2}(x)=e^{-x}-\cos(x)+\sin(x)\tag{1.7}
\end{equation*}
also satisfies Eqn. (1.2) with $\lambda=\sqrt{\pi}$. Following is an example that involves distinct $\chi_{1}(s)$ and $\chi_{2}(x)$:
\begin{equation*}
    \chi_{1}(x)=\frac{1}{4}\sin(x)+\frac{1}{4}e^{-\sqrt{3}x/2}\sin\left(\frac{x}{2}\right)+\frac{\sqrt{3}}{4}e^{-\sqrt{3}x/2}\cos\left(\frac{x}{2}\right),\tag{1.8}
\end{equation*}
\begin{equation*}
    \chi_{1}(x)=\frac{1}{4}\sin(x)-\frac{1}{2}e^{-\sqrt{3}x/2}\sin\left(\frac{x}{2}\right).\tag{1.9}
\end{equation*}
The above example satisfies Eqn. (1.2) with $\lambda=\sqrt{\pi}/2$. Ramanujan also established a technique that helps to identify and construct a combination of the above-listed functions that satisfies the functional equations (1.2). Readers can find more details on this along with some additional examples in [\cite{2}, Chp. 15].

The paper is arranged as follows: In $\S$ II we give a brief overview of Ramanujan's Master Theorem since we are going to use it frequently. The results may also be formulated without the aid of Ramanujan's Master Theorem but the reason we are using it is because the theorem makes it easy to handle the Mellin transform of exponential sums. In $\S$ III we explore some examples whose Mellin transform satisfies the following two functional equations:
\begin{equation*}
    \tilde{X}_1(s)=\frac{Q^s\sigma^2}{{X}(1-s)}\tilde{X}_2(1-s)\quad\mathrm{and}\quad \tilde{X}_1(s)=\frac{\sigma^2}{{X}(1-s)}\tilde{X}_2(1-s).\tag{1.10}
\end{equation*}
The reader can refer to Eqn. (3.1)-(3.3) for notations used in the above equations. We find that certain exponential sums weighted with characters satisfy the above equations. In $\S$ IV we consider the following more general functional equation:
\begin{equation*}
    \frac{Z_1(s)}{Z_2(k-s)}=\sigma^2.\tag{1.11}
\end{equation*}
where $Z_{1}(s)$ and $Z_2(s)$ denote the Mellin transform of functions $f_1(x)$ and $f_2(x)$. We explore these functions and as a consequence, we show that the Mellin transform of the Fourier series associated with certain Dirichlet L-functions and modular forms satisfy the above functional equation. Though we have considered specific examples in the previous section, the examples in this section are more general. We end the paper with discussions, problems, and conclusions in $\S$ V.

\begin{center}
    {\large \textbf{\textsc{\S II. RAMANUJAN'S MASTER THEOREM}}}
\end{center}
 Ramanujan's Master Theorem is an incredibly powerful tool in the theory of the Mellin transform that provides an analytic expression for the Mellin transform of analytic functions. This theorem was first communicated by S. Ramanujan in his quarterly reports that he sent to G. Hardy in England in 1913 \cite{3}. Ramanujan stated that if $f$ has an expansion of the form
\begin{equation*}
f(x)=\sum_{n=0}^{\infty}(-1)^{n} \frac{\phi(n)}{n !} x^{n}\tag{2.1} 
\end{equation*}
where $\phi(n)$ has a natural and continuous extension such that $\phi(0) \neq 0$, then for $s>0$, we have
\begin{equation*}
\int_{0}^{\infty} x^{s-1}\left(\sum_{n=0}^{\infty}(-1)^{n} \frac{\phi(n)}{n !} x^{n}\right) dx =\phi(-s)\Gamma(s).  \tag{2.2} 
\end{equation*}

Ramanujan's method for deriving his master's theorem was unconventional and his theorem had a problem with convergence of the integral. Hardy established some boundaries to the value of $\phi$ and derived a theorem that is in all respects convergent. 

Following is Hardy's version of the above theorem. Let $\varphi(z)$ be an analytic (single-valued) function, defined on a half-plane $H(\delta)=\{z \in {C}: \Re (z) \geq-\delta\}$ for some $0<\delta<1 .$ Suppose that, for some $A<\pi, \phi$ satisfies the growth condition $|\phi(v+i w)|<C e^{P v+A|w|}$ for all $z=v+i w \in H(\delta)$. Let $0<x<e^{-P}$ the growth condition shows that the series $\Phi(x)=\phi(0)-x \phi(1)+x^{2} \phi(3) \ldots$ converges. Observe that $\pi / \sin \pi s$ has poles at $s=-n$ for $n=0,1,2 \ldots$ with residue $(-1)^{n}$. The above integral converges absolutely and uniformly for $c \in(a, b)$ and $0<a<b<\delta$.The residue theorem yields
\begin{equation*}
{\Phi(x)=\frac{1}{2 \pi i} \int_{c-i \infty}^{c+i \infty} \frac{\pi}{\sin s \pi} \phi(-s) x^{-s} d s} \tag{2.3}
\end{equation*}
for any $0<c<\delta$.  Using Mellin inversion formula, $\forall 0<\Re s<\delta$, we get
\begin{equation}
\int_{0}^{\infty} x^{s-1}\left\{\phi(0)-x \phi(1)+x^{2} \phi(3) \ldots\right\} d x=\frac{\pi}{\sin s \pi} \phi(-s). \tag{2.4} 
\end{equation}
The substitution $\phi(u) \rightarrow \phi(u) / \mathrm{\Gamma}(u+1)$ in Eqn. (2.4) establishes Ramanujan's master theorem in its original form (Eqn. (2.2)). A brief history of Ramanujan's Master Theorem can be found in \cite{4}. And some analogues of Ramanujan's Master Theorem can be found in \cite{5} and \cite{6}.

 Now, replacing $x$ with $mx$ in Eqn.(2.2) and summing on $m$ from $1$ to $\infty$ gives
\begin{equation*}
\int_{0}^{\infty} x^{s-1}\sum_{m=1}^{\infty}f(mx) dx =\phi(-s)\Gamma(s)\zeta(s).  \tag{2.5} 
\end{equation*}
Furthermore, if we multiply a character $\chi(m)$ on both sides before summing, then we get
\begin{equation*}
\int_{0}^{\infty} x^{s-1}\sum_{m=1}^{\infty}\chi(m)f(mx) dx =\phi(-s)\Gamma(s)L(\chi, s).  \tag{2.6} 
\end{equation*}
We will use the above equation frequently throughout the sequel, mostly for the case where $\phi(m)=1$ for all $m$.

\begin{center}
{\large     \textsc{\textbf{\S III. FUNCTIONAL EQUATION INVOLVING CHARACTER SUMS}}}
\end{center}
We begin by defining a new set of notations. Let $\Re(s)>0$,
\begin{align*}
    &\tilde{X}_1(s)=\int_{0}^{\infty}{x^{s-1}\sum_{m=1}^{\infty}\chi(m)\xi_{1}(mx)dx},\tag{3.1}\\& \tilde{X}_2(s)=\int_{0}^{\infty}{x^{s-1}\sum_{m=1}^{\infty}\chi^\dagger(m)\xi_{2}(mx)dx},\tag{3.2}
\end{align*}
and
\begin{equation*}
    X(s)=\int_{0}^{\infty}{x^{s-1}\eta(x)dx}\tag{3.3}
\end{equation*}
where $\chi$ and $\chi^\dagger$ are some Dirichlet characters. We consider the problem of identifying and constructing the functions $\xi_1(x), \xi_2(x)$ and $\eta(x)$ such that the following functional equation is satisfied:
\begin{equation*}
    \tilde{X}_1(s)=\frac{Q^s\sigma^2}{{X}(1-s)}\tilde{X}_2(1-s).\tag{3.4}
\end{equation*}
Here, $\sigma$ and $Q$ are constants that are independent of $s$. We will later see how there's a possibility to eliminate the factor of $Q^s$ from the functional equation.

\quad\,\,\textit{{\textbf{Example 3.1.}}} 
Let
\begin{equation*}
    \left\{ {\begin{array}{*{20}{c}}
{\chi (m) = {\chi ^\dag }(m) = 1,}\\
{{\xi _1}(x) = {\xi _2}(x) = {e^{ - x}},}\\
{\eta \left( x \right) = \sin x.}
\end{array}} \right.\tag{3.5}
\end{equation*}
Notice that $\xi_{1}(x)$ and $\xi_2(s)$ has expansion of the form $(2.1)$ with $\phi(m)=1$ and $X(s)$ is, by definition, just the Mellin transform of $\sin(x)$. Therefore, using Ramanujan's Master Theorem (2.5) and Melin transform of $\sin(x)$ we get
\begin{equation*}
    \tilde{X}_{1}(s)=\tilde{X}_{2}(s)=\int_{0}^{\infty}x^{s-1}\sum_{m=1}^{\infty}e^{-mx}dx=\Gamma(s)\zeta(s)\quad\mathrm{and}\quad X(s)=\Gamma(s)\sin\left(\frac{\pi s}{2}\right).\tag{3.6}
\end{equation*}

Also, we know that the zeta function satisfies the following functional equation for $\Re(s)>0$ \cite{7}:
\begin{equation*}
    \zeta(s)=2^s\pi^{s-1}\sin\left(\frac{\pi s}{2}\right)\Gamma(1-s)\zeta(1-s).\tag{3.7}
\end{equation*}
Substituting the above values in Eqn.(3.4) and using the above functional equation gives
\begin{align*}
    \frac{\tilde{X}_{1}(s)X(1-s)}{\tilde{X}_{2}(1-s)}&=\frac{\Gamma(s)\zeta(s)\Gamma(1-s)\sin\left(\frac{\pi (1-s)}{2}\right)}{\Gamma(1-s)\zeta(1-s)}\\&={\Gamma(s)\sin\left(\frac{\pi s}{2}\right)\Gamma(1-s)\sin\left(\frac{\pi (1-s)}{2}\right)}2^s\pi^{s-1}\\&=\frac{\pi}{\sin(\pi s)}\sin\left(\frac{\pi s}{2}\right)\cos\left(\frac{\pi s}{2}\right)2^s\pi^{s-1}\\&=\left(\frac{\pi}{2}\right)2^s\pi^{s-1}\\&=\frac{(2\pi)^s}{2}\tag{3.8}
\end{align*}
which implies
\begin{equation*}
    {\tilde{X}_{1}(s)}=\frac{(2\pi)^s}{2X(1-s)}\tilde{X}_{2}(1-s).\tag{3.9}
\end{equation*}

Therefore, $Q=2\pi$ and $\sigma^2={1}/{{2}}$.

\quad\,\,\textit{\textbf{Example 3.2.}} Let $\chi(m)$ be a primitive character modulo $q$ such that $\chi(-1)=1$ and let $\chi^\dagger(m)=\bar{\chi}(m)$. Let
\begin{equation}
   \left\{ {\begin{array}{*{20}{c}}
{\chi(m)=\chi(m\,\mathrm{mod}\,q), \chi(-1)=1, \chi^{\dagger}(m)=\bar{\chi}(m)}\\{{\xi _1}(x) = {\xi _2}(x) = {e^{ - x}}}\\
{\eta \left( x \right) = \sin x}
\end{array}} \right.\tag{3.10}
\end{equation}
Notice that both $\xi_{1}(x)$ and $\xi_2(x)$ admit an expansion of the form $(2.1)$ with $\phi(m)=1$ and $X(s)$ is, by definition, just the Mellin transform of $\sin(s)$. Therefore, using Ramanujan's Master Theorem (2.6) and Melin transform of $\cos(x)$ we get
\begin{align*}
     &\tilde{X}_{1}(s)=\int_{0}^{\infty}x^{s-1}\sum_{m=1}^{\infty}\chi(m)e^{-mx}dx=\Gamma(s)L(s,\chi),\tag{3.11}\\&  \tilde{X}_{2}(s)=\int_{0}^{\infty}x^{s-1}\sum_{m=1}^{\infty}\bar{\chi}(m)e^{-mx}dx=\Gamma(s)L(s,\bar{\chi}).\tag{3.12}
\end{align*}

We know from the properties of L-functions \cite{7} that if $\chi(m)$ is a primitive character modulo q, then for all $\Re(s)>0$, we have
\begin{equation*}
    L(s,\chi)=\epsilon(\chi)2^s\pi^{s-1}q^{\frac{1}{2}-s}\Gamma(1-s)\sin\left(\frac{\pi(s+\kappa)}{2}\right)L(1-s.\bar{\chi})\tag{3.13}
\end{equation*}
where $\epsilon(\chi)=\tau(\chi)i^{-\kappa}q^{-\frac{1}{2}}$ and $\tau(\chi)$ is the Gauss sum defined by 
\begin{equation*}
    \tau(\chi)=\sum_{a=1}^{q}e^{2\pi ia/q}\chi(a).\tag{3.14}
\end{equation*}
It has the property that $|\tau(\chi)|=\sqrt{q}$, so $|\epsilon(\chi)|$=1. Furthermore,  $\kappa=0$ if $\chi(-1)=1$ and $\kappa=1$ if $\chi(-1)=-1$. 

Therefore,
\begin{align*}
    \frac{\tilde{X}_{1}(s)X(1-s)}{\tilde{X}_{2}(1-s)}&=\frac{\Gamma(s)L(s,\chi)\Gamma(1-s)\sin\left(\frac{\pi (1-s)}{2}\right)}{\Gamma(1-s)L(1-s,\bar{\chi})}\\&={\Gamma(s)\sin\left(\frac{\pi s}{2}\right)\Gamma(1-s)\sin\left(\frac{\pi (1-s)}{2}\right)}\epsilon(\chi)2^s\pi^{s-1}q^{\frac{1}{2}-s}\\&={\frac{\pi}{\sin(\pi s)}\sin\left(\frac{\pi s}{2}\right)\cos\left(\frac{\pi s}{2}\right)}\epsilon(\chi)2^s\pi^{s-1}q^{\frac{1}{2}-s}\\&=\left(\frac{\pi}{2}\right)\epsilon(\chi)2^s\pi^{s-1}q^{\frac{1}{2}-s}\\&=\frac{\tau(\chi)}{2}\left(\frac{2\pi}{q}\right)^s\quad\quad(\because \epsilon(\chi)=\tau(\chi)/\sqrt{q}\,\mathrm{when}\, \kappa=0)\tag{3.15}
\end{align*}
which implies
\begin{equation*}
     {\tilde{X}_{1}(s)}=\frac{\tau(\chi)}{2X(1-s)}\left(\frac{2\pi}{q}\right)^s\tilde{X}_{2}(1-s).\tag{3.16}
\end{equation*}

Therefore, $Q=2\pi/q$ and $\sigma^2={{\tau(\chi)}/{{2}}}$.

\quad\,\,\textit{\textbf{Example 3.3.}} In the above notation, if we let $\chi(m)$ to be a primitive character modulo q such that $\chi(-1)=-1$, then $\kappa=1$. In that case, if we set $\chi^\dagger(m)=\chi(m)$ and $\eta(x)=\cos x$ and keep everything same, then we get
\begin{equation*}
     {\tilde{X}_{1}(s)}=\frac{\tau(\chi)}{2iX(1-s)}\left(\frac{2\pi}{q}\right)^s\tilde{X}_{2}(1-s).\tag{3.17}
\end{equation*}

Therefore, $Q=2\pi/q$ and $\sigma^2={{\tau(\chi)}}/{{2i}}$.

\quad\,\,\textit{\textbf{Example 3.4.}} For $\Re(s)>0$, we define the Davenport-Heilbronn L-function \cite{8} as
\begin{equation}
   L_\mathrm{DH}(s,\sigma)=\frac{1-i\alpha}{2}L(s,\sigma)+\frac{1+i\alpha}{2}L(s,\bar{\sigma}),\,\,\mathrm{where}\,\,\alpha=\frac{\sqrt{10-2\sqrt{5}}-2}{\sqrt{5}-1}\tag{3.18} 
\end{equation}
and $\sigma$ is the character mod $5$ with $\sigma(2)=i$. The above series, when expanded, gives
\begin{equation*}
    L_\mathrm{DH}(s,\sigma)=1+\frac{\alpha}{2^s}-\frac{\alpha}{3^s}-\frac{1}{4^s}+\frac{0}{5^s}+...\tag{3.19}
\end{equation*}
It has been shown that the  Davenport-Heilbronn L-function satisfies the following functional equation:
\begin{equation*}
    L_\mathrm{DH}(s,\sigma)=5^{\frac{1}{2}-s}2^s\pi^{s-1}\Gamma(1-s)\cos\left(\frac{\pi s}{2}\right)L_\mathrm{DH}(1-s, \sigma).\tag{3.20}
\end{equation*}

Now, if
\begin{equation}
   \left\{ {\begin{array}{*{20}{c}}
   \chi(m)=\bar{\chi}(m)={(1-i\alpha)\sigma(m)}/{2}+({1+i\alpha})\bar{\sigma}(m)/{2}\\
{{\xi _1}(x) = {\xi _2}(x) = {e^{ - x}}}\\
{\eta \left( x \right) = \cos x}
\end{array}} \right.\tag{3.21}
\end{equation}
then using Ramanujan's Master Theorem (2.6) we get
\begin{align*}
    \tilde{X}_{1}(s)=\tilde{X}_{2}(s)&=\int_{0}^{\infty}x^{s-1}\sum_{m=1}^{\infty}\left(\frac{1-i\alpha}{2}\sigma(m)+\frac{1+i\alpha}{2}\bar{\sigma}(m)\right)e^{-mx}dx\\&=\Gamma(s)L_\mathrm{DH}(s,\sigma).\tag{3.22}
\end{align*}

Substituting the above values in Eqn.(3.4) and using the above functional equation gives
\begin{align*}
    \frac{\tilde{X}_{1}(s)X(1-s)}{\tilde{X}_{2}(1-s)}&=\frac{\Gamma(s)L_\mathrm{DH}(s,\sigma)\Gamma(1-s)\cos\left(\frac{\pi (1-s)}{2}\right)}{\Gamma(1-s)L_\mathrm{DH}(1-s,{\sigma})}\\&={\Gamma(s)\cos\left(\frac{\pi s}{2}\right)\Gamma(1-s)\cos\left(\frac{\pi (1-s)}{2}\right)}5^{\frac{1}{2}-s}2^s\pi^{s-1}\\&={\frac{\pi}{\sin(\pi s)}\cos\left(\frac{\pi s}{2}\right)\sin\left(\frac{\pi s}{2}\right)}5^{\frac{1}{2}-s}2^s\pi^{s-1}\\&=\left(\frac{\pi}{2}\right)5^{\frac{1}{2}-s}2^s\pi^{s-1}\\&=\frac{\sqrt{5}}{2}\left(\frac{2\pi}{5}\right)^s\tag{3.23}
\end{align*}
which implies,
\begin{equation*}
    {\tilde{X}_{1}(s)}=\frac{\sqrt{5}}{2X(1-s)}\left(\frac{2\pi}{5}\right)^s\tilde{X}_{2}(1-s).\tag{3.24}
\end{equation*}

Therefore $Q=2\pi/5$ and $\sigma^2={\sqrt{5}/2}$.

\quad\,\,\textit{\textbf{Example 3.5.}} Let $\Re(s)>0$. Consider the following two functions \cite{8}:
\begin{equation}
\sigma(s)=(1+5^{\frac{1}{2}-s})\zeta(s)=1+\frac{1}{2^s}+\frac{1}{3^s}+\frac{1}{4^s}+\frac{1+\sqrt{5}}{5^s}+...\tag{3.25}
\end{equation}
\begin{equation*}
    L(s,\chi)=1-\frac{1}{2^s}-\frac{1}{3^s}+\frac{1}{4^s}+\frac{0}{5^s}\tag{3.26}
\end{equation*}
where $\chi$ is a unique character mod $5$ with $\chi(2)=-1$. Both functions satisfy the following functional equation:
\begin{equation*}
    F(s)=5^{\frac{1}{2}-s}2^s\pi^{s-1}\sin\left(\frac{\pi s}{2}\right)\Gamma(1-s)F(1-s).\tag{3.27}
\end{equation*}

Now we consider the series $\sigma(s)$. Using Ramanujan's Master Theorem, it can be immediately observed that the following identity holds:
\begin{equation*}
    \tilde{X}_{1}(s)=\tilde{X}_{2}(s)=:\int_{0}^{\infty}x^{s-1}\left(\sum_{m=1}^{\infty}e^{-mx}+\sqrt{5}e^{-5mx}\right)dx=\Gamma(s)(1+5^{\frac{1}{2}-s})\zeta(s).\tag{3.28}
\end{equation*}
Therefore, we get
\begin{align*}
    \frac{\tilde{X}_{1}(s)X(1-s)}{\tilde{X}_{2}(1-s)}&=\frac{\Gamma(s)\sigma(s)\Gamma(1-s)\sin\left(\frac{\pi (1-s)}{2}\right)}{\Gamma(1-s)\sigma(1-s)}\\&={\Gamma(s)\sin\left(\frac{\pi s}{2}\right)\Gamma(1-s)\sin\left(\frac{\pi (1-s)}{2}\right)}5^{\frac{1}{2}-s}2^s\pi^{s-1}\\&=\frac{\pi}{\sin(\pi s)}\sin\left(\frac{\pi s}{2}\right)\cos\left(\frac{\pi s}{2}\right)5^{\frac{1}{2}-s}2^s\pi^{s-1}\\&=\left(\frac{\pi}{2}\right)5^{\frac{1}{2}-s}2^s\pi^{s-1}\\&=\frac{\sqrt{5}}{2}\left(\frac{2\pi}{5}\right)^s\tag{3.29}
\end{align*}
We would get a similar result for Eqn. (3.26).

Therefore, both of the following set of values  
\begin{equation}
   \left\{ {\begin{array}{*{20}{c}}
\chi(m)=\chi^\dagger(m)=1\\{{\xi _1}(x) = {\xi _2}(x) = {e^{ - x}+\sqrt{5}e^{-5x}}}\\
{\eta \left( x \right) = \sin x}
\end{array}} \right.\tag{3.30}
\end{equation}
and
\begin{equation}
   \left\{ {\begin{array}{*{20}{c}}
\chi(m)=\chi^\dagger(m)=\chi(\mathrm{mod}\,5), \chi(2)=-1\\{{\xi _1}(x) = {\xi _2}(x) = {e^{ - x}}}\\
{\eta \left( x \right) = \sin x}
\end{array}} \right.\tag{3.31}
\end{equation}
are solutions to the Eqn. (3.4). So get
\begin{equation*}
    {\tilde{X}_{1}(s)}=\frac{\sqrt{5}}{2X(1-s)}\left(\frac{2\pi}{5}\right)^s\tilde{X}_{2}(1-s).\tag{3.32}
\end{equation*}

\quad\,\,\textbf{\textit{Example. 3.6.}} A theorem due to H. Hamburger \cite{8} states that Riemann's zeta function is determined by the functional equation (3.7).  Hence, if we wish to produce other Dirichlet series
satisfying some functional equation, then it is necessary to change the functional equation (3.7) somehow. Take for example the Davenport-Heilbronn L-function as defined in example 3.4. The corresponding functional equation differs from Eqn. (3.7) in two ways. First, the factor of $5^{-s+\frac{1}{2}}$ has been introduced. Second, sine has been replaced with a cosine function. The latter one is unnecessary.

Now, let 
\begin{equation*}
    f(s)=\sum_{m=1}^{\infty}\frac{b_n}{n^s}\tag{3.33}
\end{equation*}
be a Dirichlet series defining a meromorphic function on the whole complex plane and let $s=\sigma+it$. If $d$ is some natural number greater than one then 
\begin{equation}
    1\pm \frac{\sqrt{d}}{d^s}=\pm d^{-s+\frac{1}{2}}\left(1\pm \frac{\sqrt{d}}{d^s}\right).\tag{3.34}
\end{equation}
This results provide us with Dirichlet polynomials satisfying the functional equation $f(s)=\pm d^{-s+\frac{1}{2}} f(1-s)$. Now let $A=a_1a_2a_3...a_r$ be the decomposition of $A$ into the product of $r$ positive integers. Define the polynomial $P(s)$ by
\begin{equation}
P(s)=\prod_{i=1}^{r}\left(1+\frac{\sqrt{a_j}}{a_{j}^{s}}\right).\tag{3.35}
\end{equation}
Then by definition, $P(s)$ satisfies the functional equation $P(s)= \epsilon A^{-s+\frac{1}{2}} P(1-s)$. The sign of each $\sqrt{a_j}$ can be either taken positive or negative. If an odd number of signs appear in $\sqrt{a_j}$ then $\epsilon=-1$ and if an even number of signs appear in $\sqrt{a_j}$ then $\epsilon=1$. Now suppose that $f(s)$ satisfies the functional equation $f(s)=\delta(s)f(1-s)$. Now define a new Dirichlet series $g(s)=P(s)f(s)$. Then $g(s)$ satisfies the following functional equation
\begin{equation}
    g(s)=\pm A^{-s+\frac{1}{2}}\delta(s)g(1-s).\tag{3.36}
\end{equation}

Notice that in example 3.1, we can rewrite the functional equation for the zeta function as $\zeta(s)=\xi(s)\zeta(1-s)$ where
\begin{equation}
    \xi(s)= 2^s\pi^{s-1}\sin\left(\frac{\pi s}{2}\right)\Gamma(1-s).\tag{3.37}
\end{equation}
It has been noted that in many cases, the quantity $\delta(s)$ is the same as $\xi(s)$ except for any change in the trigonometric function appearing in it. Thus we have two following cases: when $\delta(s)$ contains a sine function and the other when $\delta(s)$ contains a cosine function.

Using Ramanujan's Master Theorem, we have
\begin{align*}
    X_{1}(s)=X_2(s)= \int_{0}^{\infty}x^{s-1}&\sum\limits_{m = 1}^\infty  {{b_m}\prod_{i=1}^{r}({e^{ - mx}} + \sqrt {{a_i}} {e^{ - m{a_i}x}})} dx\\&=\Gamma(s)\prod_{i=1}^{r}\left(1+\frac{\sqrt{a_i}}{a_{i}^{s}}\right)f(s)=\Gamma(s)g(s).\tag{3.38}
\end{align*}

Let $\eta(x)=\sin x$. If $\delta(s)$ is of the form (3.37) then we get
\begin{align*}
    \frac{\tilde{X}_{1}(s)X(1-s)}{\tilde{X}_{2}(1-s)}&=\frac{\Gamma(s)g(s)\Gamma(1-s)\sin\left(\frac{\pi (1-s)}{2}\right)}{\Gamma(1-s)g(1-s)}\\&={\Gamma(s)\sin\left(\frac{\pi s}{2}\right)\Gamma(1-s)\sin\left(\frac{\pi (1-s)}{2}\right)}\left(\pm A^{-s+\frac{1}{2}}\right)2^s\pi^{s-1}\\&=\frac{\pi}{\sin(\pi s)}\sin\left(\frac{\pi s}{2}\right)\cos\left(\frac{\pi s}{2}\right)\left(\pm A^{-s+\frac{1}{2}}\right)2^s\pi^{s-1}\\&=\left(\frac{\pi}{2}\right)\left(\pm A^{-s+\frac{1}{2}}\right)2^s\pi^{s-1}\\&=\frac{\sqrt{\pm A}}{2}\left(\frac{2\pi}{\pm A}\right)^s.\tag{3.39}
\end{align*}

Therefore,
\begin{equation*}
    \tilde{X}_1(s)=\frac{\sqrt{\pm A}}{2X(1-s)}\left(\frac{2\pi}{\pm A}\right)^s\tilde{X}_2(1-s)\tag{3.40}
\end{equation*}
and the following is the solution:
\begin{equation*}
    \left\{ {\begin{array}{*{20}{c}}
{\chi (m) = {\chi ^\dag }(m) = b_m,}\\
{{\xi _1}(x) = {\xi _2}(x) = ({e^{ - mx}} + \sqrt {{a_1}} {e^{ - m{a_1}x}})({e^{ - mx}} + \sqrt {{a_2}} {e^{ - m{a_2}x}})...({e^{ - mx}} + \sqrt {{a_r}} {e^{ - m{a_r}x}}),}\\
{\eta \left( x \right) = \sin x.}
\end{array}} \right.\tag{3.41}
\end{equation*}

On the other hand, if $\eta(s)$ is of the form (3.37) with cosine instead of sine then we get the following set of solutions:
\begin{equation*}
    \left\{ {\begin{array}{*{20}{c}}
{\chi (m) = {\chi ^\dag }(m) = b_m,}\\
{{\xi _1}(x) = {\xi _2}(x) = ({e^{ - mx}} + \sqrt {{a_1}} {e^{ - m{a_1}x}})({e^{ - mx}} + \sqrt {{a_2}} {e^{ - m{a_2}x}})...({e^{ - mx}} + \sqrt {{a_r}} {e^{ - m{a_r}x}}),}\\
{\eta \left( x \right) = \cos x.}
\end{array}} \right.\tag{3.42}
\end{equation*}

The Mellin transform of above presented exponential sums weighted with characters are valid for $\Re(s)>0$ but the Mellin transform of $\sin x$ and $\cos x$ are together only valid in the range $0<\Re(s)<1$. Therefore, in every example presented above, the overall functional equation is only satisfied when $0<\Re(s)<1$.

\quad\,\,The functional equation (3.4) does not seem that interesting as compared to the functional equation that Ramanujan worked with, perhaps due to the factor of $Q^s$. There is, however, a way to get rid of the $Q^s$ factor easily using Ramanujan's Master Theorem. And therefore it is possible to obtain solutions for the following functional equation:
\begin{equation*}
    \tilde{X}_{1}(s)=\frac{\sigma^2}{X_{2}(1-s)}\tilde{X}_{2}(1-s).\tag{3.43}
\end{equation*}
where $0<\Re(s)<1$. Consider example 3.1. Suppose that $\xi_1(x)$ and $\xi_2(x)$ has expansion of the form
\begin{equation*}
    \xi_1(x)=\sum_{n=0}^{\infty}(-1)^{n} \frac{\phi(n)}{n !} x^{n},\quad \xi_2(x)=\sum_{n=0}^{\infty}(-1)^{n} \frac{\psi(n)}{n !} x^{n}\tag{3.44}
\end{equation*}
then a factor of $\frac{\phi(s)}{\psi(s)}$ would appear on the right-hand side of Eqn. (3.9). Letting $\phi(m)=(2\pi)^{-m}$ or $\psi(m)=(2\pi)^{m}$, i.e., letting $\xi_1(x)=e^{\frac{-x}{2\pi}}$ or $\xi_2(x)=e^{-2\pi x}$ the factor of $(2\pi)^s$ can be get ridden of. In what follows, 
\begin{equation*}
    \left\{ {\begin{array}{*{20}{c}}
{\chi (m) = {\chi ^\dag }(m) = 1,}\\
{{\xi _1}(x) = e^{-\frac{x}{2\pi}}, {\xi _2}(x) = {e^{ - x}},}\\
{\eta \left( x \right) = \sin x.}
\end{array}} \right. \quad\mathrm{and}\quad 
    \left\{ {\begin{array}{*{20}{c}}
{\chi (m) = {\chi ^\dag }(m) = 1,}\\
{{\xi _1}(x)= e^{-x}, {\xi _2}(x) = {e^{ - 2\pi x}},}\\
{\eta \left( x \right) = \sin x.}
\end{array}} \right.\tag{3.45}
\end{equation*}
are both solutions to the Eqn. (3.43).

This can be similarly done in other examples. This would give us the class of functions whose Mellin transform would satisfy the functional equation (3.43).

\begin{center}
{\large     \textsc{\textbf{\S IV. FUNCTIONAL EQUATION INVOLVING MODULAR FORMS}}}
\end{center}
What we saw in the previous section were just warm-up examples, presented to give some basic insight into the much larger panorama that we are going to discuss in this section.

 We begin by defining a new set of notations. Let $\Re(s)>0$,
\begin{align*}
    &{Z}_1(s)=\int_{0}^{\infty}{x^{s-1}f_1(x)dx},\tag{4.1}\\& {Z}_2(s)=\int_{0}^{\infty}{x^{s-1}f_2(x)dx}.\tag{4.2}
\end{align*}
Here we aim to have a more general approach. We do not assume that $f_{1}(x)$ and $f_{2}(x)$ have an expansion of a particular form as we did in the previous section.

In this section, we are going to investigate the functions $f_{1}(x)$ and $f_{2}(x)$ that satisfy the following functional equation:
\begin{equation*}
    \frac{Z_1(s)}{Z_2(k-s)}=\sigma^2\tag{4.3}
\end{equation*}
where $k$ is any real number and $\sigma$ is a parameter independent of $s$.

\quad\,\,\textit{{\textbf{Example 4.1.}}} 
Let $\left\{a_n\right\}$ and $\left\{b_n\right\}$ be two sequences of complex numbers that satisfy the condition $a-n, b_n=O(n^\epsilon)$ as $n\to \infty$ for some $\epsilon>0$. Let $\lambda>0$, $k\in \mathbb{R}$ and $\gamma\in\mathbb{C}$. For $\sigma>\epsilon+1$ let
\begin{equation*}
    \phi(s)=\sum_{n=1}^{\infty}\frac{a_n}{n^s}\quad\mathrm{and}\quad  \psi(s)=\sum_{n=1}^{\infty}\frac{b_n}{n^s}\tag{4.4}
\end{equation*}
and
\begin{equation*}
    \Phi=\left(\frac{\lambda}{2\pi}\right)^s\Gamma(s)\phi(s)\quad\mathrm{and}\quad \Psi=\left(\frac{\lambda}{2\pi}\right)^s\Gamma(s)\psi(s).\tag{4.5}
\end{equation*}

Let 
\begin{equation*}
    f(\tau)=\sum_{n=0}^{\infty}a_n e^{2\pi i n 
 \tau/\lambda}\quad\mathrm{and}\quad  g(\tau)=\sum_{n=0}^{\infty}b_n e^{2\pi i n\tau/\lambda} \tag{4.6}
\end{equation*}
where $\tau\in\mathbb{H}$ and $\mathbb{H}$ is the upper-half plane. These two series can be thought of as the Fourier series associated with the Dirichlet series (4.4). Then
$f(\tau)=\gamma(\tau/i)^{-k}g(-1/\tau)$ is equivalent to state that $\Phi(s)+a_0/s+\gamma b_0/(k-s)$ has analytic continuation to the entire
complex plane that is entire and bounded in every vertical strip. Furthermore, $\Phi(s)=\gamma \Psi(k-s)$. This result is known as Hecke's correspondence theorem \cite{9}, \cite{10}.

Using Ramanujan's Master Theorem we have
\begin{align*}
    &{Z}_1(s)=\int_{0}^{\infty}{x^{s-1}\left\{f(ix)-a_0\right\}dx}={\Phi(s)}\quad\mathrm{and}\quad {Z}_1(s)=\int_{0}^{\infty}{x^{s-1}\left\{g(ix)-b_0\right\}dx}={\Psi(s)}.\tag{4.7}
\end{align*}
Therefore, 
\begin{equation*}
    \left\{ {\begin{array}{*{20}{c}}
{{f_1}(x) = f(ix )-a_0}\\
{{f_2}(x) = g(ix )-b_0}
\end{array}} \right.\quad\mathrm{satisfies}\quad\frac{Z_1(s)}{Z_{2}(k-s)}=\sigma^2.\tag{4.8}
\end{equation*}
Hence, $\sigma^2={\gamma}$.

Following is an example due to Bochner which is a generalization of Hecke's correspondence theorem \cite{11}-\cite{15}.

\quad\,\,\textit{{\textbf{Example 4.2.}}} Let $\left\{a_n\right\}$ and $\left\{b_n\right\}$ be two sequences of complex numbers that satisfy the condition $a_n, b_n=O(n^\epsilon)$ as $n\to \infty$ for some $\epsilon>0$. Consider the following non-constant exponential series:
\begin{equation*}
    f(\tau)=\sum_{n=0}^{\infty}a_n e^{2\pi i n 
 \tau/\lambda_1}\quad\mathrm{and}\quad  g(\tau)=\sum_{n=0}^{\infty}b_n e^{2\pi i n\tau/\lambda_2}.\tag{4.9}
\end{equation*}
where $\lambda_1, \lambda_2>0$, $k\in\mathbb{R}$ and $\gamma\in\mathbb{C}$. Let $q(\tau)$ be the log-polynomial sum, which is defined as follows:
\begin{equation*}
    q(\tau)\sum_{1\leq j\leq L}\left(\frac{\tau}{i}\right)^{\alpha_j}\sum_{0\leq t\leq M(j)}\beta(j,y)\log^t\left(\frac{\tau}{i}\right)\tag{4.10}
\end{equation*}
where $L, M(j)$ are integers and $\alpha_j, \beta(j,t)$ are complex constants. Let 
\begin{equation*}
    \Phi(s)=\left(\frac{\lambda_1}{2\pi}\right)^s\Gamma(s)\phi(s)\quad\mathrm{and}\quad \Psi(s)=\left(\frac{\lambda_2}{2\pi}\right)^s\Gamma(s)\psi(s).\tag{4.11}
\end{equation*}
Then to say that $(\tau/i)^{-k}f(-1/\tau)=\gamma g(\tau)+q(\tau)$ is equivalent to saying that $\Phi(s)$ and $\Psi(s)$ have a meromorphic continuation to the entire complex plane, with a finite number of poles in $\mathbb{C}$. Furthermore, both functions satisfy the functional equation $\Phi(k-s)=\gamma \Psi(s)$. 

It immediately follows from Ramanujan's Master Theorem, or by general argument that the following would hold:
\begin{align*}
    &{Z}_1(s)=\int_{0}^{\infty}{x^{s-1}\left\{g(ix)-b_0\right\}dx}={\Psi(s)}\quad\mathrm{and}\quad {Z}_1(s)=\int_{0}^{\infty}{x^{s-1}\left\{f(ix)-a_0\right\}dx}={\Phi(s)}\tag{4.12}
\end{align*}

Therefore, 
\begin{equation*}
    \left\{ {\begin{array}{*{20}{c}}
{{f_1}(x) = g(ix )-b_0}\\
{{f_2}(x) = f(ix )-a_0}
\end{array}} \right.\quad\mathrm{satisfies}\quad\frac{Z_1(s)}{Z_{2}(k-s)}=\sigma^2.\tag{4.13}
\end{equation*}
Hence, $\sigma^2={1/\gamma}$.

The following two examples can be found in \cite{16}, pg. 108.

\quad\,\,\textit{{\textbf{Example 4.3.}}} Let $\lambda_1, \lambda_2>0$. Suppose that $\lambda_1\lambda_2=4\cos^2(\pi/q)$ where $q\geq 3$ is an odd integer. Let $\gamma  = 1\,\mathrm{if}\,k(1 - q/2)$ is even and $\gamma  = -1\,\mathrm{if}\,k(1 - q/2)$ is odd. Assume that
\begin{equation*}
    1+\left[k\left(\frac{1}{4}-\frac{q}{2}\right)+\frac{\gamma-1}{4}\right]>0.\tag{4.14}
\end{equation*}

Let $\left\{a_n\right\}$ be a complex sequence with $a_n=O(n^\epsilon)$ as $n\to\infty$. Define
\begin{equation*}
\phi(s)=\sum_{n=1}^{\infty}\frac{a_n}{n^s} \quad\mathrm{and}\quad f(\tau)=\sum_{n=0}^{\infty}a_n e^{2\pi i n 
 \tau/\lambda_1}.\tag{4.15}
\end{equation*}
Then there exists $\Phi(s)$ defined as
\begin{equation*}
\Phi(s)=\left(\frac{\lambda_1}{2\pi}\right)^s\Gamma(s)\phi(s)\tag{4.16}
\end{equation*}
such that $\Phi(s)$ has a meromorphic continuation to the entire $s$-plane with poles from the set $\left\{0,k\right\}$. The function $\Phi(s)$ satisfies the following functional equation
\begin{equation*}
    \Phi(k-s)=e^{2\pi i k/q}(\lambda_1/\lambda_2)^{k/2-s}\Phi(s).\tag{4.17}
\end{equation*}

Using Ramanujan's Master Theorem, we have
\begin{equation*}
    Z_1(s)=\int_{0}^{\infty}x^{s-1}\left\{f\left(ix\lambda_1/\lambda_2\right)-a_0\right\}dx=\Phi(s)\left(\lambda_1/\lambda_2\right)^{-s},\tag{4.18}
\end{equation*}
and
\begin{equation*}
Z_2(s)=\int_{0}^{\infty}x^{s-1}\left\{f\left(ix\right)-a_0\right\}dx=\Phi(s).\tag{4.19}
\end{equation*}

By definition, we have
\begin{align*}
    \frac{Z_1(s)}{Z_2(k-s)}&=\frac{\Phi(s)}{\Phi(k-s)}\left({\lambda_2}/{\lambda_1}\right)^s\\&=e^{-2\pi i k/q}\left({\lambda_1}/{\lambda_2}\right)^{s-k/2}\left({\lambda_2}/{\lambda_1}\right)^s \\&= e^{-2\pi i k/q}\left({\lambda_2}/{\lambda_1}\right)^{k/2}.\tag{4.20}
\end{align*}
Therefore, 
\begin{equation*}
    \left\{ {\begin{array}{*{20}{c}}
{{f_1}(x) = f\left(ix\lambda_1/\lambda_2\right)-a_0}\\
{{f_2}(x) = f\left(ix\right)-a_0}
\end{array}} \right.\quad\mathrm{satisfies}\quad\frac{Z_1(s)}{Z_{2}(k-s)}=\sigma^2.\tag{4.21}
\end{equation*}
Hence, $\sigma^2={e^{-2\pi i k/q}(\lambda_2/\lambda_1)^{k/2}}$.

Had it not been for the factor of $\lambda_1/\lambda_2$ in the integrand of Eqn. (4.18), the Eqn. (4.20) would have contained a factor of the form $Q^s$ which we are not willing to entertain in this section. Therefore, recall the trick that we used to investigate solutions to the functional equation (3.43). A similar trick can be used here to obtain another solution.  Again using Ramanujan's Master Theorem, we have
\begin{equation*}
    Z_1(s)=\int_{0}^{\infty}x^{s-1}\left\{f\left(ix\right)-a_0\right\}dx=\Phi(s)\tag{4.22}
\end{equation*}
and
\begin{equation*}
Z_2(s)=\int_{0}^{\infty}x^{s-1}\left\{f\left(ix\lambda_1/\lambda_2\right)-a_0\right\}dx=\Phi(s)\left(\lambda_1/\lambda_2\right)^{-s}.\tag{4.23}
\end{equation*}

By definition, we have
\begin{align*}
    \frac{Z_1(s)}{Z_2(k-s)}&=\frac{\Phi(s)}{\Phi(k-s)\left({\lambda_1}/{\lambda_2}\right)^{s-k}}\\&=e^{-2\pi i k/q}\frac{\left({\lambda_1}/{\lambda_2}\right)^{s-k/2}}{\left({\lambda_1}/{\lambda_2}\right)^{s-k}} \\&= e^{-2\pi i k/q}\left({\lambda_1}/{\lambda_2}\right)^{k/2}.\tag{4.24}
\end{align*}
Therefore, 
\begin{equation*}
    \left\{ {\begin{array}{*{20}{c}}
{{f_1}(x) = f\left(ix\right)-a_0}\\
{{f_2}(x) = f\left(ix\lambda_1/\lambda_2\right)-a_0}
\end{array}} \right.\quad\mathrm{satisfies}\quad\frac{Z_1(s)}{Z_{2}(k-s)}=\sigma^2.\tag{4.25}
\end{equation*}
Hence, $\sigma^2={e^{-2\pi i k/q}(\lambda_1/\lambda_2)^{k/2}}$.

\quad\,\,\textit{{\textbf{Example 4.4.}}} Let $\left\{a_n\right\}$ and $\left\{b_n\right\}$ be two sequences of complex numbers that satisfy the condition $a-n, b_n=O(n^\epsilon)$ as $n\to \infty$ for some $\epsilon>0$.  Suppose that $\lambda_1\lambda_2=4\cos^2(\pi/q)$ where $q\geq3$ is an even integer. Let $k(1-q/2)\in 2\mathbb{Z}$ and $\gamma=\pm 1$ simultaneously. Assume that 
\begin{equation*}
    1+\left[k\left(\frac{1}{4}-\frac{q}{2}\right)+\frac{\gamma-1}{4}\right]>0\tag{4.26}
\end{equation*}
where $k\geq 2q/(q-2)$. Let $\lambda>0$, $k\in \mathbb{R}$ and $\gamma\in\mathbb{C}$. For $\sigma>\epsilon+1$ let
\begin{equation*}
    \phi(s)=\sum_{n=1}^{\infty}\frac{a_n}{n^s}\quad\mathrm{and}\quad  \psi(s)=\sum_{n=1}^{\infty}\frac{b_n}{n^s}\tag{4.27}
\end{equation*}
and
\begin{equation*}
    \Phi(s)=\left(\frac{\lambda_1}{2\pi}\right)^s\Gamma(s)\phi(s)\quad\mathrm{and}\quad \Psi(s)=\left(\frac{\lambda_2}{2\pi}\right)^s\Gamma(s)\psi(s).\tag{4.28}
\end{equation*}
Let 
\begin{equation*}
    f(\tau)=\sum_{n=0}^{\infty}a_n e^{2\pi i n 
 \tau/\lambda_1}\quad\mathrm{and}\quad  g(\tau)=\sum_{n=1}^{\infty}b_n e^{2\pi i n\tau/\lambda_2}. \tag{4.29}
\end{equation*}

Then $\Psi(k-s)=\Phi(s)$.

Therefore, 
\begin{equation*}
    \left\{ {\begin{array}{*{20}{c}}
{{f_1}(x) = f\left(ix\right)-a_0}\\
{{f_2}(x) = g\left(ix\right)-b_0}
\end{array}} \right.\quad\mathrm{satisfies}\quad\frac{Z_1(s)}{Z_{2}(k-s)}=\sigma^2.\tag{4.30}
\end{equation*}
Hence, $\sigma^2=1$.

The following example can be found in \cite{18}, pg. 112.

\quad\,\,\textit{{\textbf{Example 4.5.}}}  Let $\left\{a_n\right\}$ and $\left\{b_n\right\}$ be two sequences of complex numbers that satisfy the condition $a_n, b_n=O(n^\epsilon)$ as $n\to \infty$ for some $\epsilon>0$. To each of these sequences, we associate a parameter $t\geq 0$ and define the following two series:
\begin{equation*}
   f(t)=\sum_{n=0}^{\infty}a_ne^{-\pi n t} \quad\mathrm{and}\quad g(t)=\sum_{n=0}^{\infty}b_ne^{-\pi n t}.\tag{4.31}
\end{equation*}
We also associate a Dirichlet series to these sequences:
\begin{equation*}
  L_f(s)=\sum_{n=0}^{\infty}\frac{a_n}{n^s} \quad\mathrm{and}\quad L_g(s)=\sum_{n=0}^{\infty}\frac{b_n}{n^s}.\tag{4.32}
\end{equation*}

Now suppose that for $t>0$
\begin{equation*}
   f\left(\frac{1}{t}\right)=\omega t^kg(t) \tag{4.33}
\end{equation*}
where $\omega$ is a complex number and  $k$ is a real number. Then $L_f(s)$ and $L_g(s)$ have an analytic continuation to the entire complex plane except at the poles $0$ and $k$ and satisfy the following functional equation:
\begin{equation*}
    \pi^{-s}\Gamma(s)L_f(s)=\omega\pi^{-(k-s)}\Gamma(k-s)L_f(k-s).\tag{4.34}
\end{equation*}

Using Ramanujan's Master Theorem we get
\begin{equation*}
    Z_1(s)=\int_{0}^{\infty}x^{s-1}\left\{f(x)-a_0\right\}dx=\pi^{-s}\Gamma(s)L_f(s) \tag{4.35}
\end{equation*}
and
\begin{equation*}
    Z_s(s)=\int_{0}^{\infty}x^{s-1}\left\{g(x)-b_0\right\}dx=\pi^{-s}\Gamma(s)L_g(s). \tag{4.36}
\end{equation*}

Therefore, it immediately follows that
\begin{equation*}
    \left\{ {\begin{array}{*{20}{c}}
{{f_1}(x) = f\left(x\right)-a_0}\\
{{f_2}(x) = g\left(x\right)-b_0}
\end{array}} \right.\quad\mathrm{satisfies}\quad\frac{Z_1(s)}{Z_{2}(k-s)}=\sigma^2.\tag{4.37}
\end{equation*}
Hence, $\sigma^2=\omega$.

Now we give a brief introduction to modular forms before we get down to studying examples that involve them \cite{18}. Let $f:\mathbb{H}\to\mathbb{C}$ be a holomorphic function defined by
\begin{equation*}
    f\left(\frac{az+b}{cz+d}\right)=(cz+d)^kf(z)\quad\quad\forall \left( {\begin{array}{*{20}{c}}
a&b\\
c&d
\end{array}} \right) \in SL_2(\mathbb{Z})\tag{4.38}
\end{equation*}
where $\mathrm{SL}_2(\mathbb{Z})$ is the full modular group and $k\in\mathbb{Z}$. We say that a holomorphic function $f$ is a modular form of weight $k$ for $\mathrm{SL}_2(\mathbb{Z})$ if it satisfies Eqn. (4.38) and is holomorphic at the cusp $i\infty$. 

 In the context of modular forms, we define the notion of Fourier expansion at cusp $i\infty$ to be an expansion in terms of the factor $e^{2\pi i z}$ where $z\in \mathbb{H}$. The cusp at $i\infty$ corresponds to the limit $z\to i\infty$, or equivalently $e^{2\pi i z}\to 0$. We define the Fourier series of $f(z)$ at the cusp $i\infty$ by
\begin{equation*}
    f(z)=\sum_{n=-\infty}^{\infty}a_ne^{2\pi i n z }.\tag{{4.39}}
\end{equation*}
Furthermore, a modular form is a cusp form if $a_0=0$ in its Fourier expansion. 

Now we introduce the "slash" notation. Let $\gamma\in \mathrm{GL}_{2}^{+}(\mathbb{R})$ and suppose that
$
\gamma  = \left( {\begin{array}{*{20}{c}}
a&b\\
c&d
\end{array}} \right).
$
If $z\in\mathbb{H}$ then let $j(\gamma,z)=cz+d$. For a holomorphic function $f$ of weight $k$ we define the slash notation as
\begin{equation*}
    (f|\gamma)(z)=(\mathrm{det}\gamma)^{k/2}j(\gamma,z)^{-k}f(\gamma z)\,\,\mathrm{where}\,\,\gamma z = \left( {\frac{{az + b}}{{cz + d}}} \right).\tag{4.40}
\end{equation*}

The full modular group $\mathrm{SL}_2(\mathbb{Z})$ have the following subgroups:
\begin{equation*}
    \Gamma_0 (N) = \left\{ {\left( {\begin{array}{*{20}{c}}
a&b\\
c&d
\end{array}} \right)}\in \mathrm{SL}_2(\mathbb{Z}); c\equiv 0\,\mathrm{mod}\, N \right\},
\end{equation*}
\begin{equation*}
    \Gamma_1 (N) = \left\{ {\left( {\begin{array}{*{20}{c}}
a&b\\
c&d
\end{array}} \right)}\in \mathrm{SL}_2(\mathbb{Z}); c\equiv 0\,\mathrm{mod}\, N,\quad  d\equiv 1\,\mathrm{mod}\, N \right\},
\end{equation*}
\begin{equation*}
    \Gamma (N) = \left\{ {\left( {\begin{array}{*{20}{c}}
a&b\\
c&d
\end{array}} \right) \in \mathrm{SL}_2(\mathbb{Z});\left( {\begin{array}{*{20}{c}}
a&b\\
c&d
\end{array}} \right)\equiv\left( {\begin{array}{*{20}{c}}
1&0\\
0&1
\end{array}} \right)} \mathrm{mod}\,N \right\},
\end{equation*}
and in particular $\Gamma (N) \subset {\Gamma _1}(N) \subset {\Gamma _0}(N) \subset \mathrm{SL}_2(\mathbb{Z})$.

Let $M_k=M_k(\mathrm{SL}_2(\mathbb{Z}))$ denote the $\mathbb{C}-$vector space of modular forms of weight $k$. $M_k(\mathrm{SL}_2(\mathbb{Z}))$ contains the space of cusp forms of weight $k$ which is denoted by $S_k=S_k(\mathrm{SL}_2(\mathbb{Z}))$. Furthermore, we define
\begin{equation*}
    M_k(\Gamma_0(N),\chi)=\left\{f\in M_k(\Gamma_1(N)):f|\gamma=\chi(d)f\,\forall \gamma={\left( {\begin{array}{*{20}{c}}
a&b\\
c&d
\end{array}} \right)}\in\Gamma_0(N)\right\}
\end{equation*}
which is a vector subspace of $M_k(\Gamma_1(N))$. And in particular $ M_k(\Gamma_0(N),\chi)= M_k(\Gamma_0(N))$ if $\chi$ is a trivial character.

Now we explore some examples that involve modular forms \cite{18}.

\quad\,\,\textit{{\textbf{Example 4.6.}}} Let $f\in S_k(\Gamma_0(N),\chi)$ be a cusp form and let 
\begin{equation*}
    f(z)=\sum_{n=1}^{\infty}a_ne^{2\pi i n z}\tag{4.41}
\end{equation*}
be its Fourier expansion at the cusp $i\infty$. Let
$
    \omega = \left( {\begin{array}{*{20}{c}}
0&{ - 1}\\
N&0
\end{array}} \right)
$
and put $g=f|\omega$. It can be observed that $g\in S_{k}(\Gamma_{0}(N),\chi)$. Let
\begin{equation*}
    g(z)=\sum_{n=1}^{\infty}b_ne^{2\pi i n z}\tag{4.42}
\end{equation*}
be its Fourier expansion at the cusp $i\infty$. It is obvious by definition that
\begin{equation*}
    g(z)=N^{k/2}(Nz)^{-k}f\left(-\frac{1}{Nz}\right).\tag{4.43}
\end{equation*}

With $f$ and $g$ as defined above, we associate a Dirichlet series to the Fourier coefficients to the cusp forms $f$ and $g$:
\begin{equation*}
    L_{f}(s)=\sum_{n=1}^{\infty}\frac{a_n}{n^s}\quad\mathrm{and}\quad L_{g}(s)=\sum_{n=1}^{\infty}\frac{b_n}{n^s}.\tag{4.44}
\end{equation*}
Then $L_{f}(s)$ and $L_{g}(s)$ extend to entire functions and satisfy the functional equation $\Lambda_f(s)=i^k\Lambda_g(k-s)$ where
\begin{equation*}
    \Lambda_f(s)=\left(\frac{\sqrt{N}}{2\pi}\right)^s\Gamma(s)L_f(s)\quad\mathrm{and}\quad \Lambda_g(s)=\left(\frac{\sqrt{N}}{2\pi}\right)^s\Gamma(s)L_g(s).\tag{4.45}
\end{equation*}

With the aid of Ramanujan's Master Theorem, we get
\begin{equation*}
    Z_1(s)=\int_{0}^{\infty}x^{s-1}f\left(ix/\sqrt{N}\right)dx=\left(\frac{\sqrt{N}}{2\pi}\right)^s\Gamma(s)L_f(s),\tag{4.46}
\end{equation*}
and
\begin{equation*}
    Z_2(s)=\int_{0}^{\infty}x^{s-1}g\left(ix/\sqrt{N}\right)dx=\left(\frac{\sqrt{N}}{2\pi}\right)^s\Gamma(s)L_g(s).\tag{4.47}
\end{equation*}

Therefore, it immediately follows that
\begin{equation*}
    \left\{ {\begin{array}{*{20}{c}}
{{f_1}(x) = f\left(ix/\sqrt{N}\right)}\\
{{f_2}(x) = g\left(ix/\sqrt{N}\right)}
\end{array}} \right.\quad\mathrm{satisfies}\quad\frac{Z_1(s)}{Z_{2}(k-s)}=\sigma^2.\tag{4.48}
\end{equation*}
Hence, $\sigma^2=i^k$.

\quad\,\,\textit{{\textbf{Example 4.7.}}} Let $f\in M_k(\Gamma_0(N),\chi)$ be a modular form and $g=f|\omega$ with $\omega$ as defined in above example. Following the same notation as in example 4.6, the functions $L_f(s)$ and $L_g(s)$ have analytic continuation to the entire complex plane except at the poles $0$ and $k$. Moreover
\begin{equation*}
    \Lambda_{f}(s)+\frac{a_0}{s}+\frac{b_0i^k}{k-s}\quad\mathrm{and}\quad \Lambda_{g}(s)+\frac{b_0}{s}+\frac{a_0i^k}{k-s}\tag{4.49}
\end{equation*}
and $\Lambda_f(s)=i^k\Lambda_g(k-s)$. Therefore, it immediately follows that
\begin{equation*}
    \left\{ {\begin{array}{*{20}{c}}
{{f_1}(x) = f\left(ix/\sqrt{N}\right)-a_0}\\
{{f_2}(x) = g\left(ix/\sqrt{N}\right)-b_0}
\end{array}} \right.\quad\mathrm{satisfies}\quad\frac{Z_1(s)}{Z_{2}(k-s)}=\sigma^2.\tag{4.50}
\end{equation*}
Hence, $\sigma^2=i^k$.

\quad\,\,\textit{{\textbf{Example 4.8.}}} Let $f\in S_k(\Gamma_0(q),\chi)$ where $\chi$ is a Dirichlet character modulo $q$ having conductor $q^{*}$  Suppose that $a_n$ are Fourier coefficients of $f$ at the cusp $i\infty$. Let $\psi$ be a primitive Dirichlet character modulo $r$. Consider the twisted $L$ series
\begin{equation*}
    L_{f}(s,\psi)=\sum_{n=1}^{\infty}\frac{a_n\chi(n)}{n^s}\tag{4.51}
\end{equation*}
and
\begin{equation*}
    \Lambda_{f}(s,\psi)=\left(\frac{\sqrt{N}}{2\pi}\right)^s\Gamma(s)L_f(s,\psi).\tag.{4.52}
\end{equation*}
Let $(r,q)=1$ and $N=qr^2$. Then $\Lambda_{f}(s,\psi)$ is a bounded entire function on the vertical strip and satisfies the functional equation $\Lambda_{f}(s,\psi)=i^k\omega(\psi)\Lambda_{g}(k-s,\psi)$ where $\omega(\psi)=\chi(r)\psi(q)\tau(\psi)^2/r$.

If 
\begin{equation*}
f(z)=\sum_{n=1}^{\infty}a_ne^{2\pi i n z}\tag{4.53}
\end{equation*}
then we define the twisted series by
\begin{equation*}
f_\psi(z)=\sum_{n=1}^{\infty}a_n\psi(n)e^{2\pi i n z}\tag{4.54}
\end{equation*}
which is an element of $S_{k}(\Gamma_0(N),\chi\psi^2)$ where $N$ is the east common multiple of $r, q^{*}r$ and $r^2$. 

Furthermore, let $g=f|\omega_q$ where 
\begin{align*}
    \omega_q=\left( {\begin{array}{*{20}{c}}
0&{ - 1}\\
q&0
\end{array}} \right).
\end{align*}
If $b_n$ are the coefficients of $g(z)$ in its Fourier expansion, then with the aid of Ramanujan's Master Theorem, we see that
\begin{equation*}
    Z_1(s)=\int_{0}^{\infty}x^{s-1}f\left(ix/\sqrt{N}\right)dx=\left(\frac{\sqrt{N}}{2\pi}\right)^s\Gamma(s)L_f(s,\psi),\tag{4.55}
\end{equation*}
and
\begin{equation*}
    Z_2(s)=\int_{0}^{\infty}x^{s-1}g\left(ix/\sqrt{N}\right)dx=\left(\frac{\sqrt{N}}{2\pi}\right)^s\Gamma(s)L_g(s,\psi)\tag{4.56}
\end{equation*}

Therefore, it immediately follows that
\begin{equation*}
    \left\{ {\begin{array}{*{20}{c}}
{{f_1}(x) = f_\psi\left(ix/\sqrt{N}\right)}\\
{{f_2}(x) = g_\psi\left(ix/\sqrt{N}\right)}
\end{array}} \right.\quad\mathrm{satisfies}\quad\frac{Z_1(s)}{Z_{2}(k-s)}=\sigma^2.\tag{4.48}
\end{equation*}
Hence, $\sigma^2=i^k\omega(\psi)$.

We do not assert that we have explored \textit{every} example that satisfies Eqn. (4.3). The examples provided here are merely a selection of Fourier series associated with certain Dirichlet L-functions and modular forms among the numerous instances documented in the literature. The primary objective of this section was to demonstrate that, despite certain Dirichlet L-functions and modular forms obeying distinct functional equations, the Mellin transform of their Fourier series adheres to the same functional equation.

\begin{center}
{\large     \textsc{\textbf{\S V. CONCLUSION}}}
\end{center}
In this paper, we have built upon Ramanujan's work and identified solutions for Eqns. (3.4), (3.43), and (4.3). In the process, we have shown that the Fourier series associated with certain Dirichlet L-functions and modular forms adhere to the same functional equation. Consequently, we now have solutions to four types of functional equations: the one examined by Ramanujan, Eqn. (1.2), and the three additional equations explored in this paper. It seems that only a finite number of these equations have solutions, as the existence of solutions to these functional equations appears to be governed by underlying symmetries of the functions that satisfy them.

It can be observed that in each case, we obtain a constant $\sigma^2$ that is independent of $s$ and appears to be unique in most of the examples we examined. If we were to identify a functional equation whose solutions are elements of the Selberg class, and if the constant $\sigma^2$
  proves to be unique in this context as well, it could aid in the unique classification of the elements of the Selberg class. However, a challenge in this approach lies in the fact that the functional equations of certain L-functions often involve higher powers of the gamma function or multiple gamma factors. Finding the appropriate solution in terms of the Fourier series of coefficients associated with these L-functions may require the Mellin inversion of these multiple gamma factors, a problem we believe can be addressed using the results from the paper \cite{19}.

Two significant problems remain unresolved: As discussed in $\S$ I, for Eqn. (1.2), any linear combination of known solutions yields another valid solution. This naturally raises the question of whether an appropriate combination of the Fourier series could also serve as a solution to Eqn. (4.3). Ramanujan successfully developed a theory based on Eqn. (1.3), where he established a relationship between the Mellin transform of $\phi(x)$ and $\psi(s)$ with $X_1(s)$ and $X_2(s)$, as presented in Eqn. (1.4). Moreover, he independently constructed $\phi(s)$. Consequently, it is reasonable to inquire whether a similar relationship can be formulated for Eqn. (4.3).

We further propose that by replacing $s$ with a function of $s$, say $f(s)$, in Eqns. (1.2), (3.4), (3.43), and (4.3), a broader class of solutions may be uncovered. We leave this possibility open for future research.

\textbf{Acknowledgement:} The author would like to thank Zachary Bradshaw for suggesting remarks and corrections on the final draft of this manuscript.

\end{document}